\documentclass[12pt]{amsart}

\usepackage{amsfonts}
\usepackage{amssymb,}

\catcode`\@=11

 \def\AMSTeXfeatures{\Plainheads 
   \let\current@vert=\AMS@vert}

 \def\Plainheads{\sh@ftdiam=0.05em
   \getlabeldims
   \let\vshaftfill=\plnvsolidfill
   \let\hshaftfill=\plnhsolidfill
   \let\th@rhead=\plnrhead
   \let\th@lhead=\plnlhead
   \let\th@dnhead=\plndnhead
   \let\th@uphead=\plnuphead}
 
 \def\glet{\global\let}

 \def\LaTeXfeatures{\catcode`\@=11
   \ifx\@clnwd\undefined \nol@g
      \input ltxcode.tex \dol@g \fi
   \ltxheads \let\current@vert=\new@vert
   \providelto \catcode`\@=\active}

 \def\nol@g{\def\wlog{\edef\garbage}}
 \def\dol@g{\let\wlog=\wl@g} \let\wl@g=\wlog
 \nol@g 

 \newbox\ltobox
 \def\providelto{{\setbox\z@=
   \hbox{$\to$}\minharrlen=\wd\z@
   \global\setbox\ltobox=\hbox{$\activeat>>>$}}
   \def\lto{\mathrel{\copy\ltobox}}}

 \def\ltxheads{\sh@ftdiam=\@wholewidth
   \getlabeldims
   \let\vshaftfill= \ltxvsolidfill
   \let\hshaftfill=\ltxhsolidfill
   \let\th@rhead=\ltxrhead
   \let\th@lhead=\ltxlhead
   \let\th@dnhead=\ltxdnhead
   \let\th@uphead=\ltxuphead}
 {\catcode`\@=\active
   \gdef@#1{\csname #1\string@at\endcsname}
   \glet\activeat=@}
 \def\def@#1{\expandafter\def\csname #1@at\endcsname}

 \def@>#1>#2>{\@rrow R{#1}{#2}}
 \def@<#1<#2<{\@rrow L{#1}{#2}}
 \def@ V#1V#2V{\@rrow V{#1}{#2}}
 \def@ A#1A#2A{\@rrow A{#1}{#2}}
 \def@/#1/#2/#3/{\@rrow{#1}{#2}{#3}}
 \def@.{\ifodd\row\ifmmode\noharrow
     \else\leavevmode.\spacefactor3000 \fi
   \else\novarrow\fi}
 \def@={\ifodd\row\harrow\hequalfill{}{}%
   \else\varrow\vequalfill{}{}\fi}
 \def@:#1{\ifx=#1\harrow\deffill{}{}%
   \else\leavevmode\null:#1\fi}
 \def@|{\current@vert}
  \def\AMS@vert{\varrow\vequalfill{}{}}
  \def\new@vert#1|#2|{\ifodd\row
   \let\nextarrow\vertexvarrow
   \else\let\nextarrow\varrow\fi
   \nextarrow\vshaftfill{#1}{#2}}
 \def@-{\ifmmode\let\next\hl@ne
   \else\let\next\AMSatdash \fi \next}
  \def\hl@ne#1-#2-{\harrow\hshaftfill{#1}{#2}}
  \def\AMSatdash{\let\next\relax\leavevmode
    \def\next@{\ifx\next-%
      \def\next-{\futurelet\next\nextii@}%
     \else\def\next{\hbox{-}}\fi\next}%
    \def\nextii@{\ifx\next-\def\next-{\hbox{---}}%
      \else\def\next{\hbox{--}}\fi\next}%
    \futurelet\next\next@}
 \def@(#1){\tweenarrows{#1}}
 \def@[#1]{\setsp@n#1\relax\activeat}
 \def\fiberbox{\hbox{$\vcenter{\hr@le\hbox{\vr@le
   \kern1ex\vbox{\kern1.2ex}\vr@le}\hr@le}$}}
  \def\hr@le{\hrule height \sh@ftdiam}
  \def\vr@le{\vrule width \sh@ftdiam}
 \def@+#1+#2+#3+{\ifodd\row \harrow{#1}{#2}{#3}%
   \else \varrow{#1}{#2}{#3}\fi}

 \def\Dnarrfill{\vequalfill\Dnhe@d}
 \def\Uparrfill{\Uphe@d\vequalfill}
 
 \def\ontofill{\rtarrfill\kern-0.3em 
   \th@rhead\kern 0.3em} 

 \def\rtarrfill{\hshaftfill\th@rhead}
 \def\ltarrfill{\th@lhead\hshaftfill}
 \def\dnarrfill{\vshaftfill\th@dnhead}
 \def\uparrfill{\th@uphead\vshaftfill}
 \def\hequalfill{\plnhfill=}
 \def\deffill{:\plnhfill=}
 \def\plnvextfill#1{\setbox\z@
   \hbox{\the\textfont3 #1}%
   \dimen@=\dp\z@\advance\dimen@\ht\z@
   \copy\z@ \kern-\dimen@ 
   \cleaders\copy\z@ \vfill
   \kern-\dimen@ 
   \box\z@}
 \def\plnhfill#1{$\m@th\mkern-1.5mu\mathord#1\mkern-6mu
    \cleaders\hbox{$\mkern-2mu\mathord#1\mkern-2mu$}\hfill
    \mkern-6mu\mathord#1\mkern-1.5mu$}
 \def\vequalfill{\plnvextfill{\char'167}}
 \def\plnvsolidfill{\plnvextfill{\char'077}}
 \def\plnhsolidfill{\plnhfill-}
 \def\ltxhsolidfill{\leaders\hrule height\topofshaft depth\botofshaft
   \hfill}
 \def\ltxvsolidfill{\leaders\vrule width\sh@ftdiam\vfill}
 \def\hdashfill{\hd@sh\wd@sh
   \xleaders \hbox{\wd@sh\hd@sh\wd@sh}\hfill
   \wd@sh\hd@sh}
 \def\vdashfill{\vd@sh\wd@sh
   \xleaders \vbox{\wd@sh\vd@sh\wd@sh}\vfill
   \wd@sh\vd@sh}
 \def\dashed{\ifinmeasureCD\else
    \ifodd\row\option{\let\hshaftfill=\hdashfill}%
   \else\option{\let\vshaftfill=\vdashfill}\fi\fi}

 \newdimen\CDstrutht  \newdimen\CDstrutdp
 \newdimen\CDstrutlen \CDstrutlen=\CDstrutht
   \advance\CDstrutlen by \CDstrutdp

 \def\CDstrut{\vrule
   height \ifnum\row=1 \z@\else\CDstrutht \fi
   depth \ifnum\row=\numrows \z@ \else\CDstrutdp \fi
   width\z@}

 \newdimen\CDarrsurr \CDarrsurr=0.375em
 \newdimen\CDdashlen
 \newdimen\CDvarrlen \CDvarrlen=1.5\baselineskip
 \newdimen\minharrlen 
  \setbox\z@\hbox{$\longrightarrow$} \minharrlen=\wd\z@
 \newdimen\minCDharrlen \minCDharrlen=2.5em 
\newdimen \minc@lwd
\def\findminc@lwd{\minc@lwd=2\CDarrsurr
  \advance\minc@lwd\minCDharrlen}

 \newdimen\sh@ftdiam

 \newdimen\labelsurr \labelsurr=1.25 em

\newcount\sp@ncnt \sp@ncnt=\@ne
\newcount\sp@ncnt@ \sp@ncnt@=\@ne
\newdimen\@rrwd \newdimen\@rrdp

 \def\adjustbot#1{\option{\advance\@rrdp#1\relax}}
 
\def\pushvertex#1{\global\p@shlen#1\relax
   \global\let\maybepush=\dopush}

 \newdimen\p@shlen \p@shlen=\z@

 \let\maybepush=\relax
 \def\dopush{\ifinmeasureCD 
   \advance\locdimen by -\p@shlen 
   \else\advance \@rrwd by -\p@shlen \fi 
   \global\let\maybepush=\relax \global\p@shlen=\z@\relax}

 \def\span@ne{\global\sp@ncnt=\@ne\relax}
 \def\setsp@n#1#2{\global\sp@ncnt=#1\relax
   \ifx\relax#2\relax\else\global\sp@ncnt@=#2\relax\fi}

 \def\plnrhead{\llap{$\rightarrow\mkern-1.5mu$}}
 \def\plnlhead{\rlap{$\mkern-1.5mu\leftarrow$}}

 \def\clap#1{\hbox to \z@{\hss #1\hss}}

 \def\plndnhead{\hbox{\the\textfont3 \char'171}}
 \def\plnuphead{\hbox{\the\textfont3 \char'170}}
 \def\Dnhe@d{\hbox{\the\textfont3 \char'177}}
 \def\Uphe@d{\hbox{\the\textfont3 \char'176}}

 \def\ltxrhead{\raise\@xisheight
   \llap{\smash{\@linefnt\@getrarrow(1,0)}}}
 \def\ltxlhead{\raise\@xisheight
   \rlap{\@linefnt\@getlarrow(-1,0)}}
 \def\ltxuphead{\setbox\z@=\rlap{%
   \kern\@halfwidth\@linefnt\char'66}%
   \copy\z@\kern-\ht\z@}
 \def\ltxdnhead{\setbox\z@=\rlap{%
   \kern\@halfwidth\@linefnt\char'77}%
   \ht\z@=\z@\box\z@}

 \def\wd@sh{\kern0.5\CDdashlen}
 \def\hd@sh{\vrule height\topofshaft depth\botofshaft
    width\CDdashlen}
 \def\vd@sh{\hrule height\CDdashlen
   depth\z@ width\sh@ftdiam}

\def\xylist{14{3434}13{2414}12{1723}%
  23{1413}34{1153}11{0867}43{0707}%
  32{0580}21{0414}31{0291}41{0}}
\newcount\tgtcnt@
\def\find@xyargs{\dimen@=\@rrdp
  \advance\dimen@ by \CDstrutlen
  \tgtcnt@=\dimen@ \dimen@=\@rrwd 
  \divide\dimen@ by \@m 
  \divide \tgtcnt@ by \dimen@ 
  \expandafter\testxy\xylist\relax
  \unitlength=\@xarg\@rrdp
  \divide\unitlength by\@yarg\relax}
\def\testxy#1#2#3{\ifnum\tgtcnt@>#3
    \@xarg=#1\relax \@yarg=#2\relax
    \let\next=\ignorerest
  \else\let\next\testxy\fi\next}
\def\ignorerest#1\relax{\relax}

\let\scalefactor=\@ne
\def\SWarrow{\find@xyargs\vector
  (-\@xarg,-\@yarg)\scalefactor\hskip-\wd\@linechar}
\def\NWarrow{\find@xyargs\vector
  (-\@xarg,\@yarg)\scalefactor\hskip-\wd\@linechar}
\def\NEarrow{\find@xyargs\vector
  (\@xarg,\@yarg)\scalefactor}
\def\SEarrow{\find@xyargs\vector
  (\@xarg,-\@yarg)\scalefactor}
\def\rightupline{\find@xyargs\@linelen=\scalefactor
     \unitlength\@sline}
\def\rightdownline{\find@xyargs\@yarg=-\@yarg\relax
     \@linelen=\scalefactor\unitlength\@sline}

\def\Sim{\ifodd\row\setbox\z@=\hbox{$\sim$}\dimen@=\ht\z@
 \advance\dimen@ by -\@xisheight
  \vbox{\box\z@\kern-\@xisheight\kern\dimen@}%
  \else\hbox{$\wr$}\fi}

\def\harrow#1#2#3{\inmeasureCDtrue\findminarrwd
  {#2}{#3}{\sp@ncnt\minharrlen}\inmeasureCDfalse\span@ne
  \mathrel{\hbox{\options\hplace{#1}\ulabel{#2}\dlabel{#3}}}}

\def\noharrow{\harrow\hfill{}{}}
\def\vertexvarrow#1#2#3{\findarrdp \@rrwd=\z@ \setsp@n\@ne\@ne
  \vbox to \z@{\kern-1.2\CDstrutht
  \rlap{\options\vplace{#1}\llabel{#2}\rlabel{#3}}\vss}}

\newif\ifinmeasureCD
\def\measurelabel#1{\setbox\z@
  \hbox{$\scriptstyle#1\kern\labelsurr$}%
  \ifdim\wd\z@>\@rrwd \@rrwd=\wd\z@\fi}
\def\findminarrwd#1#2#3{\@rrwd=#3\relax
   \measurelabel{#1}\measurelabel{#2}}
\def\findCDarrwd#1#2{\@rrwd=\minCDharrlen
   \measurelabel{#1}\measurelabel{#2}%
  }

\newcount\row \row=\@ne \newcount\col \col=\@ne 
 \newcount\numrows 
\numrows=\@ne
 \newcount\numcols
\newcount\arrspan \newdimen\vrtxhalfwd  \newbox\tempbox

\def\DANABUG{\advance\col by \@ne
 \@rrwd=\minCDharrlen
  \advance\@rrwd by \vrtxhalfwd
  \advance\@rrwd by \CDarrsurr
  \ifnum\col>\numcols \numcols=\col
     \newlocdimen{col\the\col}\locdimen=\@rrwd 
  \else \ifdim\@rrwd>\c@l \c@l=\@rrwd\fi\fi}

\def\drop#1\\{
  \findvrtxhalfsum\DANABUG\advance\row by 2 \measureinit}

\def\measureinit{\col=\@ne \vrtxhalfwd=-\CDarrsurr\arrspan=\@ne\@rrwd=\z@
   \setbox\tempbox=\hbox\bgroup$}
\def\measure{
  \let\harrow\measureCDarrow
  \let\CDCR=\measureCR 
   \findminc@lwd 
  \inmeasureCDtrue
  \row=\@ne \numcols=\z@ \measureinit}

\def\endmeasure{\findvrtxhalfsum\DANABUG
  \numrows=\row 
  \inmeasureCDfalse}

\def\newlocdimen#1{\advance\dimenc@unt by \@ne
  \ifnum\dimenc@unt<\insc@unt
     \else\errmessage{No room for the CD}\fi
  \dimendef\locdimen=\dimenc@unt
  \expandafter\dimendef\csname#1\endcsname=\dimenc@unt}

 \def\r@wc@l{\csname row\the\row col\the\col\endcsname}
 \def\c@l{\csname col\the\col\endcsname}

 \def\findvrtxhalfsum{$\egroup
  \newlocdimen{row\the\row col\the\col}
  \locdimen=\vrtxhalfwd 
  \vrtxhalfwd=0.5\wd\tempbox 
  \advance\vrtxhalfwd by \CDarrsurr
  \advance\locdimen by \vrtxhalfwd 
  \advance\@rrwd by \locdimen 
  \maybepush
  \divide\@rrwd by \arrspan\relax
  \ifdim\@rrwd<\minc@lwd
    \ifnum\col>\@ne \@rrwd=\minc@lwd\fi \fi
  \loop 
    \ifnum\col>\numcols \numcols=\col
       \newlocdimen{col\the\col}
       \locdimen=\@rrwd 
    \else \ifdim\@rrwd>\c@l \c@l=\@rrwd\fi \fi
   \ifnum\arrspan>\@ne
      \advance\arrspan by -1 \advance\col by \@ne
  \repeat }

 \def\measureCDarrow#1#2#3{\findvrtxhalfsum
   \arrspan=\sp@ncnt\relax\global\sp@ncnt=1\relax
   \advance\col by \@ne
   \findCDarrwd{#2}{#3}%
   \setbox\tempbox=\hbox\bgroup$}

 \newcount\dr@tn \dr@tn=\z@
 \def\locate#1:#2{\ifinmeasureCD\else
   \count@=-#1
   \multiply\count@ by 2
   \advance\count@ by #2
   \dimen@=\count@\@rrwd
   \ifnum\dr@tn=\@ne\relax \else\dimen@=-\dimen@ \fi
   \dimen@i=\@rrdp
   \ifnum\dr@tn>\z@\advance\dimen@i by \CDstrutlen \fi
   \dimen@i=\count@\dimen@i
   \count@=#2 \multiply\count@ by 2
   \divide\dimen@ by \count@
   \divide\dimen@i by \count@
   \lift\dimen@i\nudge\dimen@\fi}

\def\betweenCDrows{\advance\row by \@ne \col=\@ne
\options}

\def\hbegin{\hbox\bgroup\kern\c@l \kern-\r@wc@l$}
\def\hend{$\glet\maybepush\relax \CDstrut\egroup}
\def\vbegin{\setbox\tempbox=\hbox\bgroup$}
\def\vend{$\egroup\ht\tempbox=\z@\dp\tempbox\CDvarrlen
  \box\tempbox}
\def\setCD{\let\harrow=\setCDarrow
  \let\CDCR=\setCR 
  \row=\@ne \col=\@ne \hbegin}
\let\endsetCD=\hend 

\def\findarrwd{\@rrwd=\z@ \count@=\col \advance\count@ by\sp@ncnt
  \loop\ifnum\count@>\col \advance\count@ by -1
      \advance\@rrwd by\csname col\the\count@\endcsname\repeat}
\def\setCDarrow#1#2#3{\kern\CDarrsurr\advance\col by \@ne
  \findarrwd \advance\@rrwd by -\r@wc@l  
  \@rrdp=\z@ 
  \maybepush
  \advance\col by -\@ne \advance\col by \sp@ncnt \span@ne
  \hbox to \@rrwd{\options
   \@rrwd=\scalefactor\@rrwd\hss
   \hplace{#1}\ulabel{#2}\dlabel{#3}\hss}%
   \kern\CDarrsurr}

\newdimen\labspacei 
\newdimen\labspaceii 

\newdimen\@xisheight
  \@xisheight=\the\fontdimen22\textfont2
\newdimen\labelskip
  \labelskip=\the\fontdimen10\textfont3 
\newdimen\topofshaft
\newdimen\botofshaft
\newdimen\botofulabel
\newdimen\topofdlabel
\def\getlabeldims{
  \topofshaft=0.5\sh@ftdiam
  \botofshaft=\topofshaft
  \advance\topofshaft by \@xisheight  
  \advance\botofshaft by -\@xisheight  
  \botofulabel=\topofshaft
  \advance\botofulabel by \labelskip
  \topofdlabel=\botofshaft
  \advance\topofdlabel by \labelskip}

\def\ulabel{\ifnum\row=\@ne\let\next\ulabeli
   \else\let\next\ulabellap\fi\next}
\def\ulabeli#1{\vbox{
  \clap{\kern-\@rrwd$\scriptstyle#1$}%
  \kern\botofulabel}\maybeoffset}
\def\ulabellap#1{\vbox to \z@{\vss
  \clap{\kern-\@rrwd$\scriptstyle#1$}%
  \kern\botofulabel}\maybeoffset}
\def\dlabel{\ifnum\row=\numrows\let\next\dlabeli
   \else\let\next\dlabellap\fi\next}
\def\dlabeli#1{\vtop{\kern\topofdlabel
  \clap{\kern-\@rrwd$\scriptstyle#1$}%
  }\maybeoffset}
\def\dlabellap#1{\vbox to \z@{\kern\topofdlabel
  \clap{\kern-\@rrwd$\scriptstyle#1$}%
  \vss}\maybeoffset}
\def\rlabel#1{\vbox to \z@{\vss
  \rlap{\kern\labelskip$\scriptstyle#1$}%
  \vss\kern-\@rrdp}\maybeoffset}
\def\llabel#1{\vbox to \z@{\vss
  \llap{$\scriptstyle#1$\kern\labelskip}%
  \vss\kern-\@rrdp}\maybeoffset}
\def\swlabel#1{\vtop{\kern0.5\@rrdp
  \llap{$\scriptstyle#1$\kern\labelskip\kern-0.5\@rrwd}
  }\maybeoffset}
\def\nwlabel#1{\vbox{
  \llap{$\scriptstyle#1$\kern\labelskip\kern-0.5\@rrwd}%
  \kern-0.5\@rrdp}\maybeoffset}
\def\selabel#1{\vtop{\kern0.5\@rrdp
  \rlap{\kern0.5\@rrwd\kern\labelskip$\scriptstyle#1$}%
  }\maybeoffset}
\def\nelabel#1{\vbox{
  \rlap{\kern0.5\@rrwd\kern\labelskip$\scriptstyle#1$}%
  \kern-0.5\@rrdp}\maybeoffset}
\def\cplace#1{\vbox to \z@{\vss
  \clap{$#1$\kern-\@rrwd}%
  \kern-\@rrdp\vss}\maybeoffset}
\def\hplace#1{\hbox to \@rrwd{#1}\maybeoffset}
\def\vplace#1{\clap{\vbox to \z@{#1\kern-\@rrdp}}\maybeoffset}

\newdimen\nudgeamount \nudgeamount=\z@
\newdimen\liftamount \liftamount=\z@
\let\maybeoffset\relax
\newbox\offsetbox \newdimen\lastheight
\def\dooffset{
  \setbox\offsetbox=\lastbox \lastheight=\ht\offsetbox 
  \setbox\offsetbox=\vbox{\kern-\liftamount\box\offsetbox}%
  \ht\offsetbox=\lastheight
  \kern\nudgeamount\box\offsetbox\kern-\nudgeamount
  \global\nudgeamount=\z@ \global\liftamount=\z@
  \glet\maybeoffset=\relax}
\def\nudge#1{\ifinmeasureCD\else
  \global\advance\nudgeamount#1\relax
  \global\let\maybeoffset\dooffset\fi}
\def\lift#1{\ifinmeasureCD\else
  \global\advance\liftamount#1\relax
  \global\let\maybeoffset\dooffset\fi}

\def\findarrdp{\@rrdp=\CDvarrlen
  \ifnum\sp@ncnt@>1
    \advance\@rrdp by \CDstrutlen
    \multiply\@rrdp by \sp@ncnt@
    \advance\@rrdp by -\CDstrutlen \fi
 }

\def\varrow#1#2#3{\ifnum\sp@ncnt>\@ne 
     \sp@ncnt@=\sp@ncnt\relax\fi
  \findarrdp \@rrwd=\z@ 
  \kern\c@l
   \hbox to \z@{\options
   \@rrdp=\scalefactor\@rrdp
    \hss\vplace{#1}\llabel{#2}\rlabel{#3}\hss}%
  \global\advance\col by \@ne \setsp@n\@ne\@ne
  }

\def\novarrow{\varrow\vfill{}{}}

\def\tweenarrows#1{\findarrwd \findarrdp \setsp@n\@ne\@ne
  \rlap{\options\cplace{#1}}}

\def\usarrow #1#2#3{\dr@tn=\@ne
  \findarrwd \findarrdp \setsp@n\@ne\@ne 
  \rlap{\options\cplace{#1}\nwlabel{#2}\selabel{#3}}%
  \dr@tn=\z@}
\def\dsarrow #1#2#3{\dr@tn=\tw@
  \findarrwd \findarrdp \setsp@n\@ne\@ne 
  \rlap{\options\cplace{#1}\swlabel{#2}\nelabel{#3}}%
  \dr@tn=\z@}
 \def\@rrow#1{\csname #1@rrow\endcsname}
 \def\R@rrow{\harrow \rtarrfill}
 \def\L@rrow{\harrow \ltarrfill}
 \def\V@rrow{\varrow \dnarrfill}
 \def\A@rrow{\varrow \uparrfill}
 \def\SE@rrow{\dsarrow \SEarrow}
 \def\NW@rrow{\dsarrow \NWarrow}
 \def\SW@rrow{\usarrow \SWarrow}
 \def\NE@rrow{\usarrow \NEarrow}
 \def\DS@rrow{\dsarrow \dnslope}
 \def\US@rrow{\usarrow \upslope}
 \def\upslope{\find@xyargs
       \@linelen=\unitlength\@sline}
 \def\dnslope{\find@xyargs\@yarg=-\@yarg\relax
       \@linelen=\unitlength\@sline}

\newtoks\optionlist 
\optionlist={}
\let\options\relax
\def\dooptions{\the\optionlist\global\optionlist={}%
  \glet\options=\relax}
\def\option#1{\ifinmeasureCD\else
  \glet\options=\dooptions
  \global\optionlist=\expandafter{\the\optionlist\relax#1}\fi}
\def\wider#1{\ifinmeasureCD\else
   \option{\advance\@rrwd by #1}\fi}
\def\deeper#1{\ifinmeasureCD\else
   \option{\advance\@rrdp by #1}\fi}

{\def\\{\global\let\sptoken= }\\ }

\def\CR{\futurelet\nexttok\testCR}
\def\testCR{\ifx\nexttok\sptoken
   \let\next\eatspaceCR\else\let\next\CDCR\fi\next}
\def\eatspaceCR#1 {\CR}
\def\measureCR{\ifx\nexttok\endmeasure\let\nextCR\relax
    \else\let\nextCR\drop\fi\nextCR}
\def\setCR{\ifodd\row
  \ifx\nexttok\endsetCD\else\hend\betweenCDrows\vbegin\fi
  \else\vend\betweenCDrows\hbegin\fi}

\countdef\dimenc@unt=11
\def\CD#1\endCD{
   \begingroup\let\\=\CR
  \m@th\offinterlineskip
   \measure#1\endmeasure\null\,\vcenter{\setCD#1\endsetCD}\,
   \endgroup
    }

\ifx\@clnwd\undefined \nol@g\else\catcode`\ =14\relax\fi
 \font\@linefnt=line10 
 \newcount\@tempcnta
 \newcount\@tempcntb
 \newdimen\@tempdima
 \newdimen\@tempdimb
 \newdimen\@wholewidth
 \newdimen\@halfwidth
   \@wholewidth\fontdimen8\@linefnt \@halfwidth .5\@wholewidth
 \newdimen\unitlength
 \newcount\@xarg
 \newcount\@yarg
 \newcount\@yyarg
 \newbox\@linechar
 \newdimen\@linelen
 \newdimen\@clnwd
 \newdimen\@clnht
 \newif\if@negarg
 
 \def\@whilenoop#1{}

 \def\@whiledim#1\do #2{\ifdim #1\relax#2\@iwhiledim{#1\relax#2}\fi}

 \def\@iwhiledim#1{\ifdim #1\let\@nextwhile=\@iwhiledim 
         \else\let\@nextwhile=\@whilenoop\fi\@nextwhile{#1}}

 \def\@sline{\ifnum\@xarg< 0 \@negargtrue \@xarg -\@xarg \@yyarg -\@yarg
   \else \@negargfalse \@yyarg \@yarg \fi
 \ifnum \@yyarg >0 \@tempcnta\@yyarg \else \@tempcnta -\@yyarg \fi
 \ifnum\@tempcnta>6 \@badlinearg\@tempcnta0 \fi
 \ifnum\@xarg>6 \@badlinearg\@xarg 1 \fi
 \setbox\@linechar\hbox{\@linefnt\@getlinechar(\@xarg,\@yyarg)}%
 \ifnum \@yarg >0 \let\@upordown\raise \@clnht\z@
    \else\let\@upordown\lower \@clnht \ht\@linechar\fi
 \@clnwd=\wd\@linechar
 \if@negarg \hskip -\wd\@linechar \def\@tempa{\hskip -2\wd\@linechar}\else
      \let\@tempa\relax \fi
 \@whiledim \@clnwd <\@linelen \do
   {\@upordown\@clnht\copy\@linechar
    \@tempa
    \advance\@clnht \ht\@linechar
    \advance\@clnwd \wd\@linechar}%
 \advance\@clnht -\ht\@linechar
 \advance\@clnwd -\wd\@linechar
 \@tempdima\@linelen\advance\@tempdima -\@clnwd
 \@tempdimb\@tempdima\advance\@tempdimb -\wd\@linechar
 \if@negarg \hskip -\@tempdimb \else \hskip \@tempdimb \fi
 \multiply\@tempdima \@m
 \@tempcnta \@tempdima \@tempdima \wd\@linechar \divide\@tempcnta \@tempdima
 \@tempdima \ht\@linechar \multiply\@tempdima \@tempcnta
 \divide\@tempdima \@m
 \advance\@clnht \@tempdima
 \ifdim \@linelen <\wd\@linechar
    \hskip \wd\@linechar
   \else\@upordown\@clnht\copy\@linechar\fi}
 
 \def\@getlinechar(#1,#2){\@tempcnta#1\relax\multiply\@tempcnta 8
 \advance\@tempcnta -9 \ifnum #2>0 \advance\@tempcnta #2\relax\else
 \advance\@tempcnta -#2\relax\advance\@tempcnta 64 \fi
 \char\@tempcnta}
 
 \def\vector(#1,#2)#3{\@xarg #1\relax \@yarg #2\relax
 \@tempcnta \ifnum\@xarg<0 -\@xarg\else\@xarg\fi
 \ifnum\@tempcnta<5\relax
 \@linelen=#3\unitlength
 \ifnum\@xarg =0 \@vvector 
   \else \ifnum\@yarg =0 \@hvector \else \@svector\fi
 \fi
 \else\@badlinearg\fi}
 
 \def\@svector{\@sline
 \@tempcnta\@yarg \ifnum\@tempcnta <0 \@tempcnta=-\@tempcnta\fi
 \ifnum\@tempcnta <5
   \hskip -\wd\@linechar
   \@upordown\@clnht \hbox{\@linefnt  \if@negarg 
   \@getlarrow(\@xarg,\@yyarg) \else \@getrarrow(\@xarg,\@yyarg) \fi}%
 \else\@badlinearg\fi}
 
 \def\@getlarrow(#1,#2){\ifnum #2 =\z@ \@tempcnta='33\else
 \@tempcnta=#1\relax\multiply\@tempcnta \sixt@@n \advance\@tempcnta
 -9 \@tempcntb=#2\relax\multiply\@tempcntb \tw@
 \ifnum \@tempcntb >0 \advance\@tempcnta \@tempcntb\relax
 \else\advance\@tempcnta -\@tempcntb\advance\@tempcnta 64
 \fi\fi\char\@tempcnta}
 
 \def\@getrarrow(#1,#2){\@tempcntb=#2\relax
 \ifnum\@tempcntb < 0 \@tempcntb=-\@tempcntb\relax\fi
 \ifcase \@tempcntb\relax \@tempcnta='55 \or 
 \ifnum #1<3 \@tempcnta=#1\relax\multiply\@tempcnta
 24 \advance\@tempcnta -6 \else \ifnum #1=3 \@tempcnta=49
 \else\@tempcnta=58 \fi\fi\or 
 \ifnum #1<3 \@tempcnta=#1\relax\multiply\@tempcnta
 24 \advance\@tempcnta -3 \else \@tempcnta=51\fi\or 
 \@tempcnta=#1\relax\multiply\@tempcnta
 \sixt@@n \advance\@tempcnta -\tw@ \else
 \@tempcnta=#1\relax\multiply\@tempcnta
 \sixt@@n \advance\@tempcnta 7 \fi\ifnum #2<0 \advance\@tempcnta 64 \fi
 \char\@tempcnta}
\catcode`\ =10

\dol@g 
\catcode`\@=\active
\LaTeXfeatures

\begin{document}
\numberwithin{equation}{section}

\newtheorem{theorem}{Theorem}[section]
\newtheorem{lemma}[theorem]{Lemma}

\newtheorem*{theorema}{Theorem A}
\newtheorem*{theorema1}{Theorem A${}^\prime$}
\newtheorem*{theoremb}{Theorem B}
\newtheorem*{corc}{Corollary C}

\newtheorem{prop}[theorem]{Proposition}
\newtheorem{proposition}[theorem]{Proposition}
\newtheorem{corollary}[theorem]{Corollary}
\newtheorem{corol}[theorem]{Corollary}
\newtheorem{conj}[theorem]{Conjecture}
\newtheorem{sublemma}[theorem]{Sublemma}
\newtheorem{quest}[theorem]{Question}

\theoremstyle{definition}
\newtheorem{defn}[theorem]{Definition}
\newtheorem{example}[theorem]{Example}
\newtheorem{examples}[theorem]{Examples}
\newtheorem{remarks}[theorem]{Remarks}
\newtheorem{remark}[theorem]{Remark}
\newtheorem{algorithm}[theorem]{Algorithm}
\newtheorem{question}[theorem]{Question}
\newtheorem{problem}[theorem]{Problem}
\newtheorem{subsec}[theorem]{}
\newtheorem{clai}[theorem]{Claim}

\def\toeq{{\stackrel{\sim}{\longrightarrow}}}
\def\into{{\hookrightarrow}}

\def\kp{Val}


%

\title  [Algebraic logic and logical geometry] {Algebraic logic and logical geometry in arbitrary varieties of algebras}

\author[ Boris Plotkin] { Boris Plotkin}

\address{Boris Plotkin: Institute of
Mathematics, Hebrew University, 91904, Jerusalem, ISRAEL}
\email{plotkin@macs.biu.ac.il}



\maketitle

\long\def\symbolfootnote[#1]#2{\begingroup%
\def\thefootnote{\fnsymbol{footnote}}\footnote[#1]{#2}\endgroup}

The paper consists of two parts. The first part is devoted to logic for universal algebraic geometry. The second one deals with problems and some results. It may be regarded as a brief exposition of some ideas from the book \cite{PAP}.

\section{\bf Logic for Universal Algebraic Geometry}

\subsection{ Getting started}\label{sub:1}

For me personally, the topic of this paper originates from two main sources. The first one was my interest to mathematical models in knowledge theory, knowledge bases and databases, see, in particular,  \cite{Plotkin_Haz}, \cite{Plotkin_UA-AL-Datab}, \cite{Plotkin_AGinFOL}, \cite{Pl_Pl}.

Let us describe briefly how the bridge between algebra and knowledge theory works. We consider the following three components of knowledge:

\begin{enumerate}
\item
 A syntactical part of knowledge,  based on a language
of the given logic, is {\it the description of  knowledge.}
\item
 {\it The subject of  knowledge}  is an object
in the given applied field, i.e., an object for which we determine
knowledge. In algebraic terms the subject of knowledge is presented by an algebra $H$ in a variety $\Theta$ or by a model over this algebra.
\item
 {\it The content of  knowledge} (its semantics). Using some abuse of language we can consider the content of knowledge as a reply to the query to a knowledge base.
\end{enumerate}

A certain category of formulas in algebraic logic is related to a knowledge description. We consider a knowledge description  as a system of equations or, more generally, a system of formulas. It corresponds the knowledge content, which consists of  solutions of the given system. These solutions are presented by points in the corresponding affine space. The category of knowledge content having definable sets in the affine space as objects, is defined in a natural way.  Any passage from the knowledge description to the knowledge content is determined by a functor from the first category to the second. This functor depends on the subject of knowledge.

All these notions are defined with respect to some variety of algebras $\Theta$. In algebraic setting a variety $\Theta$ is a counter-part of  the notion of a data type defined for databases.

The described approach to knowledge theory motivates studies in logical geometry.

 Another inspiration is related to acquaintance  with the works of E.Rips and Z.Sela particularly   presented at the Amitsur Seminar in Jerusalem, and also with the works of V.Remeslennikov, O.Kharlampovich, A.Myasnikov and others. In these papers algebraic geometry in free groups has been developed. In parallel, within many years I was influenced by general viewpoints of A.Tarski and A.Maltsev on elementary theories of algebras and models.

 Both sources described above gave rise to the idea of universal algebraic geometry (UAG). In UAG we try to transit from classical algebraic geometry associated with the variety of commutative and associative algebras over a field to geometry and logic in an arbitrary variety $\Theta$ and fixed algebra $H\in \Theta$. The case when $\Theta$ is the variety of groups and $H$ is a free group in $\Theta$ remains of principal importance in UAG. We shall note that the general viewpoint from the positions of UAG to the classical variety provides some new tasks in the classical geometry as well.

\subsection{ Points and spaces of points}\label{sub:} Let an algebra $H\in \Theta$ be given. Take a finite set of variables $X=\{x_1, \ldots, x_n\}$. Define points as maps of the form $\mu: X\to H$. Each point $\mu$ determines the sequence $(h_1,\ldots, h_n)$, where $h_i=\mu(x_i)$. Since we are working in a given variety $\Theta$, one also can define a point $\mu$ as a homomorphism
$$
\mu: W(X)\to H,
$$
where $W(X)$ is the free in $\Theta$ algebra over $X$.

The set $Hom(W(X),H)$ of all homomorphisms from $W(X)$ to $H$ is regarded as an affine space or, what is the same, a space of points.

Every point defined in such a way has the classical kernel $Ker(\mu)$ and, as we will see later, the logical kernel $LKer(\mu)$.

Along with free in $\Theta$ algebras $W(X)$ we consider also algebras of formulas $\Phi(X)$ which are also associated with the given $\Theta$. We leave the precise definition of the algebra $\Phi(X)$ till Subsection \ref{sub:6}. Right now we can note that $\Phi(X)$ is an extended boolean algebra which means that $\Phi(X)$ is a boolean algebra with the operations $\exists x$, $x\in X$ called existential quantifiers and with nullary operations of the form  $w\equiv w'$, where $w, w'\in W(X)$, called equalities. There is a bunch of axioms regulating $\Phi(X)$.

We view all formulas $u\in \Phi(X)$ as equations. In particular, the formulas of the  form $w\equiv w'$, where $w, w'\in W(X)$ are equations, since they are the elements of   $\Phi(X)$.

Let a point $\mu: W(X)\to H$ be given. The logical kernel $LKer(\mu)$ consists of all
formulas $u\in \Phi(X)$ valid on the point $\mu$ (see Definition{\ref{de:lk}). This is an ultrafilter in the algebra $\Phi(X)$.


\subsection{Route map}\label{sub:3}

 We fix an infinite set of variables $X^0$ and the system $\Gamma$ of all finite subsets $X$ in $X^0$. Everywhere in the sequel we will ground on this infinite system $\Gamma$ of finite sets instead of one infinite $X^0$. This is forced by the necessity  to use logic in universal algebraic geometry. So, assume that all $X$ in  algebras $W(X)$ and $\Phi(X)$ belong to $\Gamma$.

 Let us list in the order of appearance the main structures we are going to introduce.

$\bullet$ -- $\Theta$,  the variety of algebras $\Theta$ we started with,

\medskip

$\bullet$ -- $\Theta^0$,   the category $\Theta^0$ of all free algebra $W(X)$. Morphisms in $\Theta^0$ are presented by homomorphisms $s:W(X)\to W(Y)$. Algebras $W(X)\times W(X)$ is the place, where the equations $w=w'$ live.

\medskip

$\bullet$ -- $\Theta(H)$, the category of affine spaces, see \ref{sub:4},

\medskip

$\bullet$ -- $Hal_\Theta(H)$, the category of algebras of the form $Hal_\Theta^X(H)$, see \ref{sub:4}. Algebras $Hal_\Theta^X(H)$ are the very specific algebras, built on the base of affine spaces $Hom(W(X),H)$.

\medskip

$\bullet$ -- $Hal_\Theta$, the variety of multi-sorted Halmos algebras. All algebras   $Hal_\Theta(H)=(Hal^X(H), \ X\in \Gamma)$, belong to this variety, see \ref{sub:5},

\medskip

$\bullet$ -- $\widetilde\Phi=(\Phi(X), X\in \Gamma),$ the multi-sorted free in $Hal_\Theta$ over equalities  algebra of formulas, see \ref{sub:6}. Algebras $\Phi(X)$ is the place, where the formulas live.

\medskip

$\bullet$ -- $Hal_\Theta^0$, the category of algebras of formulas $\Phi(X)$,
with morphisms
$s_*:  \Phi(X)\to \Phi(Y)$, see \ref{sub:6},

\medskip

$\bullet$ -- $ Val_{H}= (Val^{X}_{H}, X\in \Gamma)$, the homomorphism calculating values of formulas from  $\Phi(X)$ in algebras $Hal_\Theta^X(H)$,  see \ref{sub:8}.

In order to do that we shall define, first, the variety of multi-sorted Halmos algebras $Hal_\Theta$. This variety $Hal_\Theta$ corresponds to the initial variety of algebras $\Theta$. Then, we construct $Hal_\Theta^0$ as a category, corresponding to a free multi-sorted algebra $\widetilde \Phi$ in $Hal_\Theta$, built over a multi-sorted system of equalities.

We will build the morphisms $s_*$ in $Hal_\Theta^0$ in such a way, that the correspondence $s\to s_*$ determines a covariant functor from $\Theta^0\to Hal_\Theta^0$. Thus, $s_*$ has to preserve the boolean structure of the algebra of formulas $\Phi(X)$ and be correlated with quantifiers and equalities. In particular,
 $$s_{*}(w\equiv w')=(s(w)\equiv s(w')).$$
For other $u\in \Phi(X)$ the formula $s_* u=v\in \Phi(Y)$ is calculated in a more complicated way. Here, one can observe that formulas $v\in \Phi(Y)$ can contain also variables from $X\neq Y$. In fact, for $v\in \Phi(Y)$ there is no canonical form, which represent $v$ through the variables.

We will start to fulfill this program from the next Section.

\subsection{Algebras and categories $Hal_\Theta(H)$}\label{sub:4}

Algebras from the variety $Hal_\Theta$ will be the main structures in our setting.  We approach to these algebras by introducing the category $Hal_\Theta(H)$. Let us start to do that.

For each algebra $H\in \Theta$ and every finite set $X\in \Gamma$ consider the algebra $Bool(W(X),H)$. This is the boolean power algebra of $Hom(W(X),H)$ with quantifiers $\exists x$, $x\in X$ and equalities.

Define, first,  quantifiers $\exists x,$  $x \in X$ on
$Bool(W(X),H)$. Let $A$ be set from $Bool(W(X),H)$. We set $\mu \in \exists x A$ if and only if there
exists $\nu \in A$ such that $\mu(y) = \nu(y)$ for every $y \in
X$, $y\neq x$. It can be checked  that $\exists x$ defined in such
a way is, indeed, an existential quantifier.

An equality $[w\equiv w']_H$ in $Bool(W(X),H)$ is defined by
$$
\mu \in [w\equiv w']_H \Leftrightarrow (w,w')\in Ker (\mu).
$$

For some reason we denote the obtained algebra by $Hal^X_\Theta(H)$. This is an example of the extended boolean algebra in the variety $\Theta$. An algebra of formulas $\Phi(X)$ is the structure of this kind.

 We will define the algebras $\Phi(X)$ in Subsection \ref{sub:6}. In Subsection \ref{sub:6} we define the important homomorphisms

 $$
  Val_H^X: \Phi(X)\to Hal_\Theta^X(H).
 $$

Define now the category $Hal_\Theta(H)$.  Its objects are just defined algebras $Hal^X_\Theta(H)$, where $H$ is given and $X\in \Gamma$. In order to define morphisms in  $Hal_\Theta(H)$, consider first the category $\Theta(H)$ of affine spaces.

The objects of $\Theta(H)$ are affine spaces $Hom(W(X),H)$. Assign to each morphism $s: W(X)\to W(Y)$ the map $\tilde s: Hom(W(Y),H)\to Hom(W(Y),H)$ defined by the rule $\tilde s(\mu)=\mu s : W(X)\to H$, for $\mu: W(Y)\to H$. These $\tilde s$ are morphisms in $\Theta(H)$.

The correspondence $W(X)\to Hom(W(X),H)$ and $ s\to \tilde s$ defines a contravariant functor $\Theta^0\to \Theta(H)$ which determines duality.

\begin{theorem} The categories $\Theta^0$ and $\Theta(H)$ are dually isomorphic under this functor if and only if $Var(H)=\Theta$.
\end{theorem}

Morphisms $s_*^H$ in $Hal_\Theta(H)$ are defined as follows. Every
homomorphism $s: W(X)\to W(Y)$ gives rise to a Boolean
homomorphism
$$
s_{*}^H: Bool(W(X),H)\to Bool(W(Y),H),
$$
defined by the rule: for each $A\subset Hom(W(X),H)$ the point
$\mu$ belongs to $s_{*}A$ if $\tilde s(\mu)=\mu s\in A$.

The defined category $Hal_\Theta(H)$ can be treated as a multi-sorted algebra
$$
Hal_\Theta(H)=(Hal^X(H), \ X\in \Gamma),
$$
\noindent
with objects as domains and morphisms
 $$
s_{*}^H: Bool(W(X),H)\to Bool(W(Y),H),
$$
as operations.

The algebra of formulas $\widetilde \Phi=Hal_\Theta^0= (\Phi(X),\ X\in \Gamma)$ is defined in a similar way.

\subsection{Signature of algebras $Hal_\Theta(H)$}\label{sub:5}

Our next aim is to describe the signature of operations for the multi-sorted algebras $Hal_\Theta(H)$. This signature should be also multi-sorted.

 Consider an
arbitrary $W(X)$ in $\Theta$,  and take the signature
$$
L_X = \{ \vee, \wedge, \neg, \exists x, M_X \} ,\ \mbox{for all } x
\in X.
$$
  Here
$M_X$ is the  set of all formulas $w\equiv w'$, $w,w'\in W(X)$
over the algebra $W(X)$.  We treat these formulas as symbols of relations of equality over $W(X)$, that is there is a map $\equiv: W(X)\times W(X)\to M_X$ which satisfies axioms of an equational predicate on $W(X)$.  These are the only symbols of relations in use. Symbols  $w\equiv w'$ can be regarded also as symbols of nullary operations.

Signature $L_X$ is the signature of the one-sorted extended boolean algebras.  Now we define the signature $L_\Theta$ for the multi-sorted algebras $Hal_\Theta(H)$.

Along with the set of symbols of equalities $M_X$ consider the set $S_{X,Y}$ of symbols of operations $s_*$ 
of the type $\tau=(X;Y)$, where $X,Y\in\Gamma$.
Symbols $s_*$ are just symbols of operations but we keep in mind that each
homomorphism $s:W(X)\to W(Y)$ induces the operation $s_*$ in $Hal_\Theta(H)$ of the
type $\tau=(X;Y)$.

By the same reason we assume that given $s: W(X)\to W(Y)$ and $s': W(Y)\to W(Z)$, the axiom
$$
(ss')_*=s_*s'_*
$$
\noindent
holds. Here the operation $(ss')_*$ has the type $\tau=(X;Z)$.

Define the signature
 $$
 L_\Theta= \{L_X, S_{X,Y}; X, Y \in \Gamma \}
$$




The signature $L_\Theta$ is multi-sorted. We consider the constructed multi-sorted algebras  $Hal_\Theta(H)$ in this signature with the natural realization of all operations from $L_\Theta$. We will take $L_\Theta$ also as the signature of an arbitrary algebra from the variety of multi-sorted algebras $Hal_\Theta$.

Now we define algebras which belong to the variety $Hal_\Theta$.



\begin{defn}\label{def:ha}
We call an algebra $\frak L = (\frak L _X, X \in \Gamma)$ in the
signature $L_\Theta$ a Halmos algebra, if

\begin{enumerate}

\item
 Every domain $\frak L _X$ is an extended Boolean algebra in the signature $L_X$.
\medskip
\item
Every mapping $s_*: \frak L _X \to \frak L _Y$ is
  a homomorphism of Boolean
algebras. Let $s: W(X) \to W(Y)$, $s': W(Y) \to W(Z)$, and let
$u\in \frak L _X$. Then $s'_*(s_*(u))=(s's)_*(u)$.

\medskip

\item The identities controlling the interaction of $s_{*}$ with
quantifiers are as follows:

\begin{enumerate}

\medskip
\item[(a)]
 $s_{1*} \exists x a = s_ {2*} \exists x
a, \ a \in \frak L (X)$, if $s_1( y) = s_2(y)$ for every $y\neq x$, $x$, $y\in X$.


 \medskip

 \item[(b)]
   $s_{*}\exists x a = \exists (s(x)) (s_*a),$ $\ a \in \frak L (X)$,
if $ s(x) = y $ and $y $ is a variable which does not belong to
the support
 of $s(x')$,
 for every
$x' \in X$, and $x' \neq x$.

\noindent
This condition means that $y$  does not participate in the shortest expression of the
element $s(x')\in W(Y)$. 

\end{enumerate}

\item
   The identities controlling the interaction of $s_{*}$ with equalities are as follows:
  \begin{enumerate}

  \medskip

  \item[(a)]
   $s_{*}(w\equiv w')=(s(w)\equiv s(w'))$.

   \medskip

  \item[(b)]
   $(s^{x}_{w})_{*}a \wedge (w\equiv w')\le (s^{x}_{w'})_{*}a$,
   where $a\in \frak L (X)$, and $s^{x}_{w}\in End(W(X))$ is defined by $s^{x}_{w}( x)=w$,
   and $s^{x}_{w}(x')=x',$ for $x'\ne x$.

 \end{enumerate}

\end{enumerate}

\end{defn}

One should not be upset with the looking complicated axioms from the items 3--4. First of all, we have already algebras $Hal_\Theta(H)$ as an example of Halmos algebras, So, one can verify in the very straightforward way that $Hal_\Theta(H)$  satisfy these axioms.  Second, the purely logical explanations on the language of first order formulas of the similar axioms for the one-sorted polyadic algebras are contained in \cite{Halmos}, see also \cite{Plotkin_UA-AL-Datab}.

\begin{defn} The variety $Hal_\Theta$ consists of the multi-sorted algebras in the signature $L_\Theta$ subject to axioms from Definition \ref{def:ha}.
\end{defn}

The conditions specified in Definition \ref{def:ha} have the form of identities and they actually  define a variety.

\subsection{ Free algebras in $Hal_\Theta$}\label{sub:6}

We construct the free in $Hal_\Theta$ algebra  $\widetilde \Phi$ in an explicit way.
Denote by $M=(M_X, X \in \Gamma)$ the multi-sorted  set of
equalities with the components $M_X$.


 Let us build the absolutely free algebra of formulas in the signature $L_\Theta^0$. Each formula in this algebra has two parameters: the length and the sort. Then we define formulas by induction.

 Each equality $w\equiv w'$ is a formula of the length
zero, and of the sort  $X$ if $w\equiv w'\in M_{X}$. Let $u$ be a
formula of the length $n$ and the sort $X$. Then the formulas
$\neg u$ and $\exists x u$ are the formulas of the same sort $X$
and  the length $(n+1)$. Further, for the given $s:W(X)\to W(Y)$
we have the formula $s_{*}u$ with the length $(n+1)$ and the sort
$Y$. Let now $u_1$ and $u_2$ be formulas of the same sort $X$ and
the length $n_1$ and $n_2$ accordingly. Then the formulas $u_1\vee
u_2$ and $u_1\wedge u_2$ have the length $(n_1+n_2+1)$ and the
sort $X$. In such a way, by induction, we define lengths and sorts
of arbitrary formulas.

We construct a big set formulas $
\mathfrak L ^{0}$.

Let $\mathfrak L ^{0}_{X}$ be the set of all formulas of the sort
$X$. Each $\mathfrak L ^{0}_{X}$ is an algebra in the signature
$L_{X}$ and we have the algebra
$$
\mathfrak L ^{0}=(\mathfrak L ^{0}_{X}, X\in \Gamma)
$$
 in the signature $L_{\Theta}$. By construction,
the algebra $\mathfrak L ^{0}$ is the absolutely free algebra of
formulas over equalities $M=(M_X, X \in \Gamma)$ concerned
with the variety of algebras $\Theta$. Its elements are considered as pure formulas.

Denote by  $ \tilde \pi$ the verbal congruence in $\mathfrak L ^{0}$
generated by the identities of Halmos algebras from Definition
\ref{def:ha}.


 Define the Halmos algebra of
formulas as
$$
\widetilde \Phi=\mathfrak L ^{0}/ \tilde \pi.
$$
It can be  written as $\widetilde \Phi=(\Phi(X), X\in \Gamma)$,
where
$$
\Phi(X)=\mathfrak L ^{0}_{X}/ \tilde \pi_{X},
$$
where each $\Phi(X)$ is an extended Boolean algebra of the sort
$X$ in the signature $L_X$.

 The algebra $\widetilde\Phi$ is the
free algebra in the variety $Hal_\Theta$  of all  multi-sorted Halmos algebras
associated with the variety of algebras $\Theta$, with the set of
free generators $M=(M_X, X\in \Gamma)$.

This approach to the free in $Hal_\Theta$ algebra through the factorization of the absolutely free algebra of formulas by the verbal congruence can be viewed as syntactic.


We could built $\widetilde\Phi$ also semantically. In each $\mathfrak L ^{0}_{X}$ the formulas of the sort $X$ are collected. Recall that for $s:W(X)\to W(Y)$ and $u\in \mathfrak L ^{0}_{X}$ the formula $v=s_*u$ lies in $\mathfrak L ^{0}_{X}$. All these formulas can be treated as pure formulas of a logic which possesses some axioms and rules of inference which correlate with the definition of the variety $Hal_\Theta$.

One can show 
 that if we factor out
component-wisely
 the algebra $\mathfrak L ^{0}$ by the
many-sorted
 Lindenbaum-Tarski congruence, then we get the same
algebra $\widetilde \Phi$. This observation provides a bridge
between syntactical and semantical description of the free
multi-sorted Halmos algebras.

Using the correspondence between multi-sorted algebras and categories we are able to define the category of algebras of formulas $Hal_\Theta^0$. Its objects are algebras $\Phi(X)$, its morphisms have the form $s_*: \Phi(X)\to \Phi(Y)$, where $s:W(X)\to W(Y)$ is a morphism in $\Theta^0$. According to Definition \ref{def:ha} the morphisms $s_*$  preserve the boolean structure of $\Phi(X)$ and correlated with quantifiers and equalities. Hence the correspondence
$$
W(X)\to \Phi(X) \quad  {\rm and} \quad  s\to s_*
$$
determines a covariant functor $\alpha: \Theta^0 \to Hal_\Theta^0$.

\subsection{ Variety $Hal_\Theta$ and algebras $Hal_\Theta(H)$}\label{sub:7}

As we have seen, the algebras $Hal_\Theta(H)$ belong to $Hal_\Theta$.  Moreover,

\begin{theorem}\label{thm:gen} Algebras $Hal_\Theta(H)$ generate the variety $Hal_\Theta$.
\end{theorem}

This means that we could define the algebra $\widetilde\Phi$ using algebras $Hal_\Theta(H)$.
We have a unique homomorphism of the algebra $\mathfrak L$ to  $Hal_\Theta(H)$. The kernel of this homomorphism is the system of identities of  $Hal_\Theta(H)$. Intersection of all kernels through all $H\in \Theta$ and coincide with the verbal congruence with respect to the variety $Hal_\Theta$

 One can prove that these algebras are simple with respect to congruences  and all simple algebras are exhausted by algebras $Hal_\Theta(H)$ and their subalgebras.

\subsection{The Value homomorphism}\label{sub:8}

The free algebra $\widetilde\Phi=(\Phi(X), X\in \Gamma)$ and an arbitrary algebra $Hal_\Theta(H)=(Hal_\Theta^X(H), \ X\in \Gamma)$ belong to the same variety $Hal_\Theta$.


We  define a homomorphism
$$
 Val_{H}: \widetilde\Phi \to Hal_\Theta(H),
 $$

\noindent
which    induces the homomorphism
$$
 Val^{X}_{H}: \Phi(X) \to Hal_\Theta^X(H)
 $$
of the one-sorted extended Boolean algebras, for every $X$ in $\Gamma$. 

Since equalities $M=(M_X, X\in\Gamma)$ freely generate   the free
multi-sorted Halmos algebra $\widetilde\Phi$, it is enough to assign an equality in $Hal_\Theta^X(H)$ to the corresponding equality $w\equiv w'$ in $\Phi(X)$.

Recall that we have defined equalities $[w\equiv w']_H$ in $Hal_\Theta^X(H)$ by
$$
\mu \in [w\equiv w']_H \Leftrightarrow (w,w')\in Ker (\mu).
$$

Hence, the element $  [w\equiv w']_H$ in $Hal_\Theta^X(H)$ is assigned to the element $w\equiv w'$ in $\Phi(X)$. This correspondence gives rise to the homomorphism of
multi-sorted Halmos algebras
$$
 Val_{H}: \widetilde\Phi \to Hal_\Theta(H).
 $$
Since $\widetilde\Phi =(\Phi(X), X\in \Gamma)$, where each
component $\Phi(X)$ is an extended Boolean algebra, the
homomorphism   $
 Val_{H}$ induces homomorphisms
$$
 Val^{X}_{H}: \Phi(X) \to Hal_\Theta^X(H),
 $$
of the one-sorted extended Boolean algebras. In particular,
$$
Val^{X}_{H}(w\equiv w')=\{ \mu \mid \mu (w)=\mu(w') \}.
$$

\begin{defn}
A point $\mu: W(X)\to H $  satisfies the formula $u\in \Phi(X)$ if  $Val^{X}_{H}(u)$ contains $\mu$.
\end{defn}

This definition has the same meaning as the standard model theoretic one.
  Now we are in a position to define formally the logical kernel of a point.

  \begin{defn}\label{de:lk}
A formula $u\in \Phi(X)$ belongs to the logical kernel $LKer(\mu)$ of a point $\mu$ if and
only if $\mu\in Val^{X}_{H}(u)$.
\end{defn}

  If $u\in \Phi(X)$ then $Val^{X}_{H}(u)$ is the set of points $\mu$ satisfying the formula $u$. This means that
  $$
  u\in LKer(\mu) \Leftrightarrow \mu\in Val^X_H(u).
    $$

 Since we consider each formula $u\in \Phi(X)$ as an "equation" and
$Val^{X}_{H}(u)$ as a value of the formula $u$ in the algebra
$Bool(W(X),H)$, then  $Val^{X}_{H}(u)$ is a set of points
$\mu:W(X)\to H$ satisfying the "equation" $u$. We call
$Val^{X}_{H}(u)$  {\it solutions of the equation} $u$. We also say
that the  formula $u$ holds true in the algebra $H$ at the point
$\mu$.

It can be verified that the logical kernel $LKer(\mu)$ is always a
boolean ultrafilter of $\Phi(X)$.

 Note that
 $$
 Ker(\mu)=LKer(\mu)\cap M_X.
 $$

Now let $Th(H)=(Th^X(H), X\in \Gamma)$ be the multi-sorted
representation of the elementary theory of $H$. We call its
component $Th^X(H)$ the {\it $X$-theory of the algebra $H$}. Since $Ker(Val_H)$ is the set of formulas satisfied by all points of $H$, we have
have:
$$
 Ker(Val_H)=Th (H),
 $$
  $$
 Ker(Val^X_H)=Th^X(H).
 $$

This means, in particular, that the algebra $\Phi(X)$ can be represented, modulo elementary theory, in the more transparent algebra $Hal_\Theta^X(H)$. We can also present the $X$-theory of the algebra $H$ as:

$$
Th^X(H)=\bigcap_{\mu: W(X)\to H} LKer(\mu).
$$

\begin{defn} An algebra $H$ is called saturated if for every finite $X$ every ultrafilter $T$ in $\Phi(X)$, which contains  $Th^X(H)$ coincides with a $Lker(\mu)$ for some $\mu:W(X)\to H$.
\end{defn}

This notion stimulates a lot of problems.

 The key diagram which relates logic of different sorts in multi-sorted case is as follows:

$$
\CD
\Phi(X) @> s_\ast >> \Phi(Y)\\
@V  \kp^X_H  VV @VV \kp^Y_H V\\
Hal_\Theta^X(H) @>s_\ast>> Hal_\Theta^Y(H)
\endCD
$$

Here the upper arrow represent the syntactical transitions in the
category $Hal_\Theta$, the lower level does the same with the
respect to semantics in $Hal_\Theta$, and the correlation is
provided by the vertical value homomorphism.

With this diagram we finish exposition of the necessary ideas from algebraic logic and switch  to a logically-geometric stuff. We defined formally the multi-sorted algebra of formulas $\tilde\Phi$ and its domains $\Phi(X)$, where $X\in \Gamma$. These domains could be informally treated as dynamic algebras of formulas, which means that formulas (elements) from a particular $\Phi(X)$ interact with formulas from other $\Phi(Y)$.

Now we have prepared all necessary information for the next Part.

\section{\bf Some results and problems}

\subsection{ Types and isotypeness}\label{sec:tp}

In model theory for each set $X$, $X=\{x_1,\ldots,x_n\}$ the notion of $X$-type ($n$-type) is defined (see \cite{Marker}).
 Given an algebra $H$ with the elementary theory $Th(H)$, a set $P$ of formulas $u$, such that all free variables in $u$ belongs to $X=\{x_1,\dots , x_n\}$ is an $X$-type ($n$-type),   if $P\cup Th(H)$
is satisfiable. Denote by $T_p^H(\mu)$ the $X$-type of the point $\mu: W(X)\to H$, i.e., the set of all formulas $u$ valid on $\mu$.  The types of the form $T_p^H(\mu)$ will be called {\it $MT$-types or model-theoretic types} of points.

We consider also {\it $LG$-types or logically-geometric types of the points $\mu: W(X)\to H$}. According to Definition 6.2 from \cite{PlAlPl} an $LG$-type of a point $\mu$ is the logical kernel $LKer(\mu)$ of the point $\mu$ in the algebra $\Phi(X)$. We will return to the geometric nature of this definition later. The following result connects $MT$-types and $LG$-types, and plays a key role in all considerations:

\begin{theorem}[\cite{Zhitom_types}]\label{th:zh}

 Let the points $\mu: W(X)\to H_1$ and $\nu: W(X)\to H_2$ be given.  The equality
$$
T_p^{H_1}(\mu)=T_p^{H_2}(\nu)
$$
holds if and only if we have
$$
LKer(\mu)=LKer(\nu).
$$
\end{theorem}

\medskip

 Now we will describe the idea of isotypeness of algebras.

\begin{defn}\label{df:iso} Two algebras $H_1$ and $H_2$ are called isotypic  if for every $X$ and every point  $\mu: W(X)\to H_1$ there exists a point $\nu: W(X)\to H_2$ such that the types of $\mu$ and $\nu$ coincide, and  for for every point  $\nu: W(X)\to H_2$ there exists a point $\mu: W(X)\to H_1$ such that the types of $\mu$ and $\nu$ coincide.
\end{defn}

$\bullet$ -- In view of Theorem \ref{th:zh} two algebras  $H_1$ and $H_2$  are isotypic with respect to $MT$-types (
$MT$-isotypic) if and only if they are  isotypic with respect to $LG$-types (
$LG$-isotypic). Thus,  in Definition \ref{df:iso} one can equally rely on $MT$-types and $LG$-types.

\medskip

$\bullet$ -- Since in the definition of isotypeness one can grounds on coincidence of logical kernels of the points, then this definition has a geometric nature and extends the notion of geometrically equivalent algebras  ( see  next Section \ref{sec:alg} and especially Theorem \ref{th:lg}  for details). One can say that Theorem \ref{th:zh} and Theorem \ref{th:lg} visualize the bridge between geometrical  and logical ideas.

\medskip

$\bullet$ -- According to Definition \ref{df:iso} isotypeness of algebras implies their elementary equivalence. So, this notion is more strong than the notion of elementary equivalence and expresses a logical property of algebras which should in many cases be closed to isomorphism.


Indeed, it is easy to see that if algebras $H_1$ and $H_2$ are isotypic, then they are locally isomorphic. This means that every finitely generated subalgebra in $H_1$ is isomorphic to a finitely generated subalgebra in $H_2$ and vice versa. Of course, locally isomorphic groups are not necessarily isomorphic and isotypic. For example every two free groups $F_m$ and $F_n$ are locally isomorphic, but not isomorphic provided $m\neq n$. As will see  $F_m$ and $F_n$  are not also isotypic.
This follows from the result of   C.~Perin and R.~Sklinos \cite{PerinSklinos} on logical homogenity (see Section \ref{sec:reg}) of a free group.


Moreover, it easy to see that if $H_1=F_n$ is a finitely generated  free group and $H_2$ is a finitely generated group, then their isotypeness implies isomorphism. This follows from the local isomorphism of $H_1$ and $H_2$, which means that $H_2$ can be viewed as a subgroup of $H_1$. Since every subgroup of $F_n$ is free, it remains to use the above mentioned  result on isotypeness of free finitely generated groups.

 We come up with the following conjecture:

\begin{conj}\label{pr:free}
Let $F_n$ be a free group of the  rank $n>1$ and $H$ be a  group. If $F_n$ and $H$ are isotypic then they are isomorphic.
\end{conj}

Recently, R.~Sklinos \cite{Sklinos_1} gave a positive answer to this problem.  His proof is also based  on the logical homogenity of a free group in the form of the following theorem of Pillay:

\begin{theorem}[Pillay]
Let $F_n$ be the free group with free generators $e_1, \ldots,  e_n$. Consider the points $\mu: W(X)\to F_n$ and $\nu: W(X)\to F_n$ defined by $\mu(x_i)=e_i$ and $\nu(x_i)=a_i$, respectively,  where $i=1,\ldots n$, and $a_i$ are arbitrary elements in $F_n$. Suppose that
$$
T_p^H(\mu)=T_p^H(\nu).
$$
Then $a_1,\ldots, a_n$ are the free generators of $F_n$.
 \end{theorem}

We see that all groups isotypic to a free finitely generated group $F_n$ should be isomorphic to $F_n$. The next problems are related to the general case:
\begin{problem}\label{pr:free2}
Let  the groups $H_1$ and $H_2$  be isotypic and $H_1$ be  finitely generated. Is it true that $H_2$ is finitely generated?
\end{problem}

\begin{problem}\label{pr:free1}
Let  $H_1$ and $H_2$ be two finitely generated isotypic groups. Are they isomorphic?
\end{problem}

Since all the theory grounds on an arbitrary variety of algebras Problems \ref{pr:free1} and \ref{pr:free2} which are formulated for the variety of all groups make sense for an arbitrary variety of algebras $\Theta$. Their solution heavily depends on the choice of $\Theta$.  So, let  $H_1$ and $H_2$ be two algebras from a variety $\Theta$.

\begin{problem}\label{pr:free2'}

Let  the algebras $H_1$ and $H_2$  be isotypic and $H_1$ be  finitely generated. Is it true that $H_2$ is finitely generated?
\end{problem}

\begin{problem}\label{pr:free1'}
Let  $H_1$ and $H_2$ be two finitely generated isotypic algebras. Are they isomorphic?

\end{problem}

Let us point out one more closely related question
\begin{problem}\label{pr:free1'}
Let  $H_1$ and $H_2$ be two finitely generated isotypic groups. What can be said about isotypeness of the groups algebras $KH_1$ and $KH_2$ where $K$ is a field.
\end{problem}

As we have mentioned above a local isomorphism does not imply isotypeness. Here is an example when local isomorphism of algebras implies isotipicity:

\begin{example} Any two  infinite dimensional vector spaces $H_1$ and $H_2$ are locally isomorphic. It can checked that their  local isomorphism implies isotypeness of $H_1$ and $H_2$.
 Take now   two non-isomorphic infinitely dimensional vector spaces. Then, in view of above they provide an
example of isotypic, locally isomorphic but not isomorphic algebras.
\end{example}

\medskip

\subsection{ Algebraic sets}\label{sec:alg}

 Recall some basic notions from equational algebraic geometry in the given variety $\Theta$. Fix a finite set $X=\{x_1,\ldots, x_n\}$ and let  $T$ be a system of equations of the form $w\equiv w'$,
$w,w'\in W(X)$.

Define
$$
T'_H = A = \{\mu : W(X) \to H \ | \
$$
$$
T \subset Ker(\mu)\}
$$

Subsets  $A$ in $Hom(W(X),H)$ of the form $A=T'_H$ are called {\it algebraic sets}.

Let $A$ be an arbitrary set of points $\mu:W(X)\to H$ in $Hom(W(X),H)$. Define
$$
A'_{H}=T=\{(w\equiv w')  \ |\ $$
$$ (w, w')\in \bigcap_{\mu\in A}
Ker(\mu)\}.
$$

The congruences $T$ of the form $T=A'_{H}$ are called $H$-closed congruences in $W(X)$. We obtain a Galois correspondence between different $T$ and $A$. If $A$ and $T$ are arbitrary sets of points and sets of equations, respectively, then $A''_H$ and $T''_H$ are the corresponding Galois closures.

Define now an important notion:

\begin{defn} Algebras $H_1$ and $H_2$ are called geometrically equivalent (or $AG$-equivalent) if for every finite $X$ and every $T$ in $W(X)$ we have
$$
T''_{H_1}=T''_{H_2}.
$$
\end{defn}

\noindent
A criterion for the geometrical equivalence is a follows :

\noindent
Consider  infinitary quasiidentities of the form $ \ (\ast) $

$$
w_1\equiv w_1'\wedge\ldots\wedge w_n\equiv w_n'\wedge\ldots $$

$$\Longrightarrow w_0\equiv w_0'.
$$

Then, $H_1$ and $H_2$ are $AG$-equivalent if and only if every quasiidentity of the form $(\ast) $ which holds in $H_1$ holds also in $H_2$ and vice versa. In particular, if $H_1$ and $H_2$ are $AG$-equivalent, then they have the same finitary quasiidentities. The converse is not true \cite{GS}, \cite{MR}.

One can define in a very natural way the category of algebraic sets over the given algebra $H$ in the variety $\Theta$.  Denote it by $K_\Theta(H)$.

It is easy to prove that if algebras $H_1$ and $H_2$ are $AG$-equivalent, then the categories $K_\Theta(H_1)$ and $K_\Theta(H_2)$
are isomorphic. In the paper \cite{Pl-St} there are necessary and sufficient conditions for specific algebras  $H_1$ and $H_2$  that provide isomorphism of the categories $K_\Theta(H_1)$ and $K_\Theta(H_2)$.

\medskip

\subsection{ Elementary (definable) sets}\label{sec:alg}

Now we switch to logical geometry. The main correspondences in logical geometry repeat that of algebraic geometry with replacing kernels of points $Ker(\mu)$ with their logical kernels $LKer(\mu)$.

If $T$ is set of formulas in the algebra $\Phi(X)$ then define:

$$
T^{L}_{H}=A=\{ \mu :W(X)\to H \mid $$
$$
T\subset LKer(\mu) \},
$$
Subsets  $A$ in $Hom(W(X),H)$ of the form $A=T^{L}_H$ are called {\it elementary or definable sets}.

Define

$$
A^{L}_{H}=T=\bigcap_{\mu\in A} LKer(\mu)
$$

\noindent
The filters $T$ of the form $T=A^{L}_{H}$ are called $H$-closed filters in $\Phi(X)$. We obtain a Galois correspondence between $T$ and $A$.

We have also that a formula $u\in \Phi(X)$ belongs to $T=A^L_H$ if and only if $A\subset Val_H^X(u)$.

\begin{defn} Algebras $H_1$ and $H_2$ are called logically equivalent (or $LG$-equivalent) if for every finite $X$ and every $T\subset \Phi(X)$ we have
$$
T^{LL}_{H_1}=T^{LL}_{H_2}.
$$
\end{defn}

Using infinitary formulas of the form $ \ (\ast) $:

$$
\bigwedge_{u\in T}u \Longrightarrow v, \quad u,v\in \Phi(X)
$$

we can formulate the criterion:

 Algebras $H_1$ and $H_2$ are $LG$-equivalent if and only if every formula of the form $(\ast) $ which holds in $H_1$ holds also in $H_2$ and vice versa. 

The main point is the following theorem which connects the model theoretic notion of isotypic algebras with the geometric notion of $LG$-equivalent algebras:

\begin{theorem}\label{th:lg}
Algebras $H_1$ and $H_2$ are isotypic if and only if they are $LG$-equivalent.
\end{theorem}

Consider  the category of elementary sets over the given algebra $H$ in the variety $\Theta$.  Denote it by $LK_\Theta(H)$.

\begin{theorem}
 If algebras $H_1$ and $H_2$ are isotypic, then the categories $LK_\Theta(H_1)$ and $LK_\Theta(H_2)$
are isomorphic.
\end{theorem}

In the situation of logical geometry we do not have yet  necessary and sufficient conditions on algebras  $H_1$ and $H_2$  that provide isomorphism of the categories $LK_\Theta(H_1)$ and $LK_\Theta(H_2)$.

\medskip

One can prove that every elementary set $A$ in $Hom(W(X),H)$ is invariant under the action of the group $Aut(H)$ in this affine space.

\subsection{ Logically regular varieties}\label{sec:reg}

Now,

 \begin{defn}\label{df:sep} Algebra $H\in \Theta$ is called logically separable (in $\Theta$) if every $H'\in \Theta$ which is not isomorphic to $H$ is not isotypic to $H$.
\end{defn}
This means that algebra $H\in \Theta$ is  logically separable (in $\Theta$) if every $H'\in \Theta$ isotypic  to $H$ is isomorphic to $H$.
Definition \ref{df:sep} says that a logically separable algebra $H$ can be distinguished by means of the logic of types. We are looking for the cases, when a free in $\Theta$ algebra $W(X)$, $X\in \Gamma$ is separable in $\Theta$.

\begin{defn} A variety $\Theta$  is called $LG$-regular (logically regular) if every free in $\Theta$ algebra $W(X)$, $X\in \Gamma$ is logically separable.
\end{defn}

It is checked that the varieties of semigroups, of inverse semigroups are $LG$-regular (see  \cite{Zhitom_types}), the  variety of abelian groups is $LG$-regular (see \cite{Zhitom_types} and \cite{Sklinos}).

More recently it was established that:
\medskip

$\bullet$ -- variety of all groups is $LG$-regular \cite{Sklinos_1},

\medskip

$\bullet$ -- variety of all nilpotent of class $c$ groups is $LG$-regular \cite{Zhitom_types}.

\medskip



\begin{problem}

Is it true that the following  varieties are  $LG$-regular:

\medskip

$\bullet$ -- Variety of solvable  groups of the derived length $d$.

\medskip

$\bullet$ -- Variety of  commutative associative algebras over a field.

\medskip

$\bullet$ -- Variety of  associative algebras over a field.

\medskip

$\bullet$ -- Variety of Lie algebras.

\medskip

$\bullet$ -- Variety of class $c$  nilpotent associative algebras.

\medskip

$\bullet$ -- Variety of class $c$  nilpotent Lie algebras.
\end{problem}

We shall specify also the following:

\begin{problem}
Find a variety $\Theta$ which is not $LG$-regular. Construct a non-regular variety of groups.
\end{problem}


\subsection{ Logically homogeneous algebras}\label{sec:homogeneous}

\begin{defn}
An algebra $H\in \Theta$ is called logically homogeneous if for every two points $\mu: W(X)\to H$ and $\nu: W(X)\to H$ the coincidence of their types implies that there exists an automorphism $\sigma$ of $H$ such that $\sigma(\mu)=\nu$.
\end{defn}

It can be seen that if an algebra $H$ is logically homogeneous, then for every $X$ every orbit of the action of $Aut(H)$ in $Hom(W(X),H$ is an elementary set for some $T$ in $\Phi(X)$.

Take an arbitrary point $\mu: W(X)\to H$, and consider $T=LKer(\mu)$ in $\Phi(X)$.  Then the elementary set $T^L_H$ defined by the ultrafilter $T$ is an orbit of $Aut(H)$ containing the point $\mu$, provided $H$ is logically homogeneous.

\begin{defn} We call a variety $\Theta$ logically perfect if every free in $\Theta$ algebra $W(X)$ is logically homogeneous.
\end{defn}

\begin{problem}\label{prob:21}
What are the varieties $\Theta$ such that an arbitrary free in $\Theta$ algebra $W(X)$, $X\in \Gamma$  is logically homogeneous.
\end{problem}

Problem \ref{prob:21} has positive solution for the variety of all groups \cite{PerinSklinos},
 for the varieties of
abelian groups and nilpotent of class $c$ groups
(\cite{Zhitom_types}, for abelian groups see also \cite{Sklinos}). It is also proved that  the torsion free hyperbolic groups are logically homogeneous \cite{Houcine}.

\begin{problem}\label{prob:23}
Is it true that the variety of  solvable  groups of the derived length $d>1$
is  logically perfect?
\end{problem}

In particular,

\begin{problem}\label{prob:24}
Is it true that the variety of metabelian  groups
is  logically perfect?
\end{problem}

Note, that up to now we don't have examples of non-logically homogeneous varieties of groups $\Theta$.

\begin{defn} Algebra $H$ is called algebraically homogeneous if for two points $\mu: W(X)\to H$, and $\nu: W(X)\to H$ the equality $Ker(\mu)=Ker(\nu)$ implies that there exists an automorphism $\sigma$ of $H$ such that $\sigma(\mu)=\nu$.
\end{defn}

It is clear that the algebraic homogenity means that every isomorphism $\alpha: A\to B $ of their finitely generated subalgebras $A$ and $B$ is realized by some $\tau \in Aut(H)$ . By definition, algebraic homogenity implies logical homogenity. Converse statement is not true. For example, the free abelian group is logically homogeneous, but no algebraically homogeneous.

 Algebraically homogeneous solvable groups  and finite groups are classified, see \cite{CF}, and \cite{CF1}, respectively.  These groups are very closed to the $QE$-groups, i.e., the groups whose elementary theory admits quantifier elimination. The similar theory is elaborated for algebraically homogeneous rings (see \cite{CSW} and references therein).




\medskip

\subsection{ Logical noetherianity}\label{sec:noeth}

 Consider now logically noetherian algebras.

\begin{defn} An algebra $H\in \Theta$ is called logically noetherian (LG-noetherian) if for every $X\in \Gamma$ and every $T\subset \Phi(X)$ there exists a finite subset $T_0\subset T$ such that

$$
T^L_H=T^L_{0H}.
$$

\end{defn}

In concern with this definition an algebra $H\in \Theta$ is called automorphically finitary if there is only a finite
number of $Aut(H)$ -orbits in $Hom(W(X),H)$ for every $X$.  If $H$ is automorphically finitary, then $H$ is $LG$-noetherian.

\begin{problem}

1. Describe automorphically finitary abelian groups.

2. Consider non-commutative automorphically finitary groups.

\end{problem}









If for the algebra $H$ there are only a finite
number of $Aut(H)$-orbits in $Hom(W(X),H)$ for every $X$, then
there are only finite number of realizable $LG$-types  in
$\Phi(X)$.

\subsection{  Addendum}

This paper is devoted to the theory, which gave rise to a system of  new notions. It arises  in a natural way and stimulates a lot of new problems in algebra and universal algebraic geometry which are concerned with logic and algebraic logic. Recent results of G.Zhitomirskii and R.Sklinos are of this kind. We also look at new papers of Z.Sela \cite{Sela} , and O.Kharlampovich and A.Myasnokov \cite{KhM_6} through the prism of the described theory.

Note the following general

\begin{theorem} If in the variety of algebras $\Theta$  is logically perfect, then $\Theta$ is $LG$-regular.
\end{theorem}

 We know that  the varieties of all groups, the variety of nilpotent of class $c$ groups, the variety of abelian groups are $LG$-perfect. Hence, these varieties are $LG$-regular.  This makes the problem if the varieties of solvable groups of the derived length $d$, of metabelian groups and of commutative associative rings with unity are perfect,  especially interesting.
 We have no examples of varieties which are not logically perfect and are not logically regular. So, it is quite important to check the varieties specified above.  More precisely, is it possible in these cases that the group of automorphisms does not act transitively on Galois closures of points?




We will finish with the problems related to elementary equivalence of algebras. These problems, however, appeared as a result of studies around isotypeness of algebras.

 \begin{problem} Find a variety of groups $\Theta$, which is distinct from the variety of all groups, such that all free in $\Theta$ finitely generated algebras are elementary equivalent, but if they are isotypic, then they are isomorphic. In particular,

 $\bullet$ -- Is the Burnside variety $B_n$ of groups satisfying the identity $x^n=1$, of this kind?

 $\bullet$ -- Is the Engel variety $E_n$ of groups satisfying the identity
 $$[x,y] \ldots,y]=1,$$
 where the commutator is taken $n$-times, of this kind?

 \end{problem}
  \begin{problem} Find two elementary equivalent groups $H_1$ and $H_2$ such that their group algebras are elementary equivalent, but not isotypic.
  \end{problem}

\end{document}